\documentclass[12pt]{article}
\usepackage{graphicx}
\usepackage{caption}
\usepackage{subcaption}
\usepackage{algorithmic}
\usepackage{algorithm}
\usepackage{amsmath}
\usepackage{amssymb}
\usepackage{amsfonts}
\usepackage{amsthm}
\usepackage{fullpage}
\usepackage{url}
\usepackage{color}
\usepackage{float}
\usepackage{dsfont}
\usepackage{setspace}
\usepackage{bm}
\usepackage{wrapfig}
\usepackage[font=small,labelfont=bf]{caption}
\usepackage{cancel}

\allowdisplaybreaks

\let\oldmarginpar\marginpar
\renewcommand\marginpar[1]{\-\oldmarginpar[\raggedleft\footnotesize #1]%
{\raggedright\footnotesize #1}}

\def\d{{\textrm{d}}}

\newtheorem{definition}{Definition}

\newcommand{\diag}{{\rm diag}}

\newcommand{\Z}{\mathbb{Z}}

\renewcommand{\d}{\,{\rm d}}

\def\I{{\mathcal I}}

\title{Time-frequency representation of nonstationary signals: the IMFogram}
	\author{Philippe Barbe\thanks{LDC and CNRS (on leave) ({\tt philippe.barbe@math.cnrs.fr})}
		\and
	Antonio Cicone\thanks{DISIM, Universit\`a degli Studi dell'Aquila, L'Aquila, ITALY. ({\tt antonio.cicone@univaq.it}) AC is a member of the Italian ``Gruppo Nazionale di Calcolo Scientifico'' (GNCS) of the Istituto Nazionale di Alta Matematica ``Francesco Severi'' (INdAM). His work was partially supported through the CSES-Limadou project of the Istituto di Astrofisica e Planetologia Spaziali (IAPS) of the Istituto Nazionale di Astrofisica (INAF).}
		\and
		Wing Suet Li\thanks{School of Mathematics, Georgia Institute of Technology, Atlanta, GA, USA. ({\tt li@math.gatech.edu})  WSL is partially supported by a Simon Foundation Collaboration Grant for Mathematicians. }
		\and
		Haomin Zhou\thanks{School of Mathematics, Georgia Institute of Technology, Atlanta, GA, USA. ({\tt hmzhou@math.gatech.edu}) HZ is supported in part by NSF under grants DMS-1830225, and ONR N00014-18-1-2852.}
	}
	\date{}

\begin{document}

	\maketitle
	
\centerline{\it Dedicated to the memory of Professor Ciprian Foia\c s}
\vskip 0.5 in

	\begin{abstract}
	    \noindent Iterative filtering methods were introduced around 2010 to improve  definitions and measurements of   structural features in signal processing. Like many applied techniques, they present considerable challenges for mathematicians to theorize their effectiveness and limitations in commercial and scientific usages.  In this paper we recast  iterative filtering methods in a mathematical abstraction more conducive to their understanding and applications. We also introduce a new visualization of simultaneous local frequencies and amplitudes. By combining a theoretical and practical exposition, we hope to stimulate efforts to understand better these methods. Our approach acknowledges the influence of Ciprian Foia\c s,   who was passionate about pure, applied,  and applications of mathematics.
	\end{abstract}	

	
\section{Introduction }\label{sec:Introduction}

A common way to analyze a complex signal is to
decompose it as a superposition of simpler ones. Many methods have been  devised,  from Fourier series to the more modern wavelets. These methods are extremely successful in dealing with stationary signals. However they are less adapted to non-stationary ones because
the decomposition is done in the frequency domain, and when transformed back to the time or spatial domain, yields  stationary components.
In contrast, the iterative filtering approach discussed in this paper offers a compromise between separating frequencies and retaining non-stationarity.

\subsection{ Overview.}

The traditional approach to time-frequency analysis is rooted in Fourier analysis and was later developed into wavelets and other frame decompositions
\cite{cohen1995time, daubechies1992ten,  flandrin1998time, ghobadi2020comparative,  mallat2012group,  anden2019joint,  daubechies1996nonlinear,  daubechies2011synchrosqueezed,
auger2013time, daubechies2016conceft}. These techniques, overviewed in \cite{wu2020current}, rely on predetermined bases or frames that are not data-driven.
Consequently, some results may depend on the chosen basis and may not reflect intrinsic properties of the signal \cite{flandrin1998time}.

These methods often assume that the signal is stationary, or near-stationary. Localization with distribution kernel models tries to remove the stationarity assumption. It attempts to balance bias, forcing to localize more, with stability, forcing to localize less and average more.

These traditional decompositions are linear, that is why the corresponding time-frequency representation are defined of linear-type.
In the last decades some method have been proposed to improve the linear-type time-frequency representations, such as the reassignment method \cite{auger2013time} and the synchrosqueezing transform \cite{daubechies2011synchrosqueezed}.

Bilinear-type techniques have been also proposed, like the Wigner-Ville distribution and the Cohen's or affine class. However these techniques proved to be limited in producing a clean and sharp time-frequency representation \cite{flandrin1998time,wu2020current}.

For all these reasons,
more data-driven decompositions have been designed. In particular, the empirical mode decomposition \cite{huang1998empirical}, which divides the signal into several simple oscillatory components. The components are obtained inductively, removing from the signal the average between its upper and lower envelopes, these envelopes being defined by interpolation between the local extremum (maximum for the upper envelope, minimum for the lower envelope). Each component, called an intrinsic mode in \cite{huang1998empirical}, may then be analyzed separately in the time-frequency domain
\cite{huang2014introduction, huang2009instantaneous}.
Informally, intrinsic modes try to equate the number of extrema and zero-crossing, making them nicely oscillating, and the average between their upper and lower envelopes is near $0$, making them well balanced.

Decomposition in empirical modes is successful in a wide range of applications
\cite{cummings2004travelling, parey2006dynamic, yu2008forecasting, huang2008review, tary2014spectral, stallone2020new}.
However, the algorithm may be unstable and may not separate well nearby frequencies, creating the so-called mode mixing problem \cite{Rilling2008frequencies}.
  As a consequence, several variants have been proposed
\cite{huang2009EEMD, yeh2010complementary, torres2011complete, rehman2009Filterbank}.
These methods received considerable attentions from the scientific community, as indicated by Huang and collaborators' papers receiving so far more than 30,000 citations on Scopus.


Intrinsic modes, produced using the aforementioned techniques, allow to overcome artificial spectrum spread caused by sudden changes in local frequencies, formalized in the Heisenberg-Gabor uncertainty principle \cite{flandrin1998time}. Regrading the actual computation of the frequency and amplitude content of each intrinsic mode, different approaches have been proposed in the literature, like in \cite{cicone2016adaptive,huang2009instantaneous,sandoval2018instantaneous}.

Hybrid methods optimizing a decomposition on a fixed basis or frame with data-driven amplitude or phase modulation have also been proposed
\cite{hou2011adaptive, hou2009variant, yu2018geometric, coifman2017carrier, gilles2013empirical, dragomiretskiy2013variational, selesnick2011resonance, meignen2007new, pustelnik2012multicomponent}.

\subsection{Paper outline. }

The work of \cite{cicone2016adaptive} highlights the vast literature on the subject, and from there one can easily get a better view of the progress that has been made during the last 20 years.  In this paper we consider the iterative filtering method and its variation using the fast Fourier transform, the {\it fast iterative filtering}, FIF algorithm, that are inspired by the empirical mode decomposition \cite{huang1998empirical}.

Twenty years ago,
Ciprian Foia\c s,  Mike Jolly, and Wing Suet Li \cite{foias} published in the area of numerical analysis. Ciprian joked that the method was communist:
one step forward and two steps backward! Here we are making a capitalist approach: after two steps forward (iterative filtering and its fast version), we take one step backward, for a broader view of the methodology!

The basic idea is to convolve the signal with a filter determined by features of the signal, and to decompose the signal into finitely many components, called {\it intrinsic mode functions} (IMFs). Since the filter is determined by the signal itself, each IMF keeps a significant amount of non-stationarity  and local frequency characteristics from the signal. This makes them easier to interpret by practitioners. Our approach is not definitive: we have a theoretical framework,   some interesting examples, a new visualization,  but complete mathematical justifications remain  elusive.

It is easy to explain the basic principle of iterative filtering method for continuous signals. However, in practice, one will most likely encounter  discrete finite time series. We may also deal with spatial signals, like air pressure on earth, moisture in a field, possibly discretized. Some signals, such as those relating gyroscopic data and other sensors may be indexed by the sphere.
The unifying framework of groups allows us to develop a method that covers many applications. Some of our developments  may be carried out in homogeneous spaces, a setting that includes the most intriguing examples, but key ideas become hidden in a more complex formalism.

We present the  theoretical framework in section 2. Section 3 specializes it to time series. A new visualization scheme, the IMFogram (pronounced like ``infogram''), is introduced in section 4.
We present some numerical examples in section 5.

\section{The theoretical framework }\label{sec:Theorectical}

Signals are viewed as real valued functions defined on a group such as $\mathbb{R}$ or $\mathbb{Z}$ for time series, $\mathbb{R}^2$ or $\mathbb{Z}^2$ for spatial signals, the circle or some discretization of it for periodic ones.  Because we need convolutions and Fourier transforms, we restrict  $G$ to be a locally compact abelian group written additively, with a Haar measure, denoted by $\lambda$. If $G$ is finite, $\lambda$ is the counting measure up to normalization. The following recalls basic abstract Fourier analysis, see \cite{loomis}.

The convolution of two functions $u$ and $v$ defined on $G$ is
\[
u\star v(g) = \int_G u(g-h) v(h) \d \lambda(h) \, ,\quad g\in G.
\]

The characters of $G$ are all the continuous homomorphisms from $G$ to the complex unit circle. The set of all characters,  $\widehat{G}$,  is an abelian group under the pointwise multiplication, called the dual group of $G$.

The Fourier transform ${\mathcal F}$ maps linearly a function $u$ defined on $G$ to the function defined on the dual group by
\[
{\mathcal F}u(\chi) = \int u(g)\overline\chi(g)\d \lambda(g),  \quad \chi\in\widehat{G}\, .
\]
It satisfies ${\mathcal F}(u\star v)  = {\mathcal F}u\,{\mathcal F}v$.

It is possible to find a Haar measure on $\widehat G$, $\d \chi$, such that the inverse of the Fourier transform is given by
\[
  {\mathcal F}^{-1}u(g) = \int_{\widehat G} u\bigl(\chi(g)\bigr) \d\chi \, , \qquad g\in G\,.
\]

We continue our discussion by defining the central player of the FIF method,  the filter. Recall that a function $w$ on $G$ is even if $w(x) = w(-x)$ for every $x$ in $G$.

\bigskip

\begin{definition}\label{def: filter}

Let $G$ be a locally compact abelian group, written additively, with Haar measure $\lambda$.
\begin{itemize}
	\item [(i)] A function $w$ on $G$ is  a filter if it is nonnegative,  even,  bounded, and $\displaystyle{\int w \d \lambda = 1}$.
	\item [(ii)] A double convolution filter $w$ is the self-convolution of a filter  $\tilde w$, that is, $w = \tilde w \star \tilde w$.
    \item [(iv)]  The size, or the length,  of a filter $w$ is the Haar measure of its support, $\ell(w) = \lambda\{\, w>0\,\} $.

\end{itemize}
\end{definition}

This definition implies that
convolutions of filters are filters.
The range of the Fourier transform of a filter is in  $[\,-1,1\,]$, and that of a double convolution filter is in $[\,0,1\,]$.

Consider a possibly non-stationary signal $s$ in $L^2(G)$. Its moving average with respect to a filter $w$ is defined as the convolution
\begin{equation}\nonumber
    \mathcal{C}_w s = w \star s \,.
\end{equation}
We subtract the moving average from the signal and obtain the variation of the signal around its $w$-moving average,
\begin{equation}\nonumber
    \mathcal{V}_ws = s - \mathcal{C}_ws.
\end{equation}
Iterating $p$-times the linear operator $\mathcal{V}_w$ we obtain the linear IMF operator,
\begin{equation}\nonumber
    \mathcal I_{w,p}=\mathcal{V}_w^p .
\end{equation}
We will discuss how to choose the filter $w$ and  this $p$ later.
Since filters are bounded,  IMF  operators are endomorphisms on both $L_1(G)$ and $L_2(G)$.

We can rewrite these operators on the dual group using the Fourier transform. Representing the pointwise multiplication of functions by a dot ($\cdot$),
\begin{equation}\nonumber
    \mathcal{F}\,\mathcal{I}_{w,p}s = (1 - \mathcal{F}w)^p \cdot \mathcal{F}s
\end{equation}

Taking the inverse Fourier transform,
\begin{equation}\nonumber
   \mathcal{I}_{w,p}s = \mathcal{F}^{-1}\bigl((1 - \mathcal{F}w)^p \cdot \mathcal{F}s\bigr) \, .
\end{equation}
The construction of IMFs is iterative. Starting with a signal $s$, we obtain the first IMF as some $\I_{w_1,p_1}s$. We then consider the remainder, $s_1 = s - \I_{w_1,p_1}s$ and construct a second IMF as some $\I_{w_2,p_2}s_1$. More generally, given a reminder $s_{n-1}$, we define the $n$-th IMF as  $\I_{w_n,p_n}s_{n-1}$ for some $w_n$ and $p_n$, and the next reminder as $s_{n} = s_{n-1} - \I_{w_n,p_n}s_{n-1}$.

Whenever $G$ and its dual are finite, all the linear operators involved so far may be written as finite matrices. We agree to list the elements of $G$ and its character group $\widehat G$ in some chosen specific order. A function $f$ on $G$ is then a vector $\bigl(f(g)\bigr)_{g\in G}$ and the Fourier transform $\mathcal F$ is the matrix ${\mathcal F} = \bigl(\overline\chi(g)\bigr)_{\chi\in \widehat G,\,g\in G}$, with its rows indexed by the character group $\widehat G$ and its columns by the group $G$. Then
\begin{equation}\nonumber
  \mathcal{V}_{w} = \mathcal{F}^{-1} \diag(1 - \mathcal{F}w) \mathcal{F} \, ,
\end{equation}

\begin{equation}\label{eq:keyEq}
  \mathcal{I}_{w,p} = \mathcal{V}_{w}^p
= \mathcal{F}^{-1}\diag(1 - \mathcal{F}w)^p \mathcal{F} \, .
\end{equation}
The $n$-th IMF is obtained by applying a linear operator of the following form to $s$:
\begin{equation}\nonumber
  \mathcal{F}^{-1} \diag(1 - \mathcal{F}w_n)^{p_n} \prod_{k=1}^{n-1}  \bigl( 1 -  \diag(1 - \mathcal{F}w_k)^{p_k} \bigr)    \mathcal{F}\, .
\end{equation}

However, the $n$-th IMF is not a linear function of $s$ since the $w_k$ and $p_k$ are determined by $s$.

Assume that $w$ is a double convolution filter,  ensuring that its Fourier transform is in $[\,0,1\,]$.  Then
 \begin{equation}\nonumber
   \lim_{p\to\infty} \diag(1 - \mathcal{F}w)^p = \diag\bigl(\mathds{1}\{\mathcal{F}w(\chi)=0\}\bigr)_{\chi\in \widehat{G}} \, .
\end{equation}
 This limit is the orthogonal projection on the subspace of functions in $L^2(\widehat{G})$ that vanish on the support of $\mathcal{F}w$. It is a
 bandpass filter when $G$ is $\mathbb{Z}_{N}$, and it preserves no non-stationarity information from the original signal.

In practice, the choice of $p$ reflects the desire that the IMFs focus on a section of  the  spectrum, that is $ \diag(1 - \mathcal{F}w)^p $ is fairly concentrated, yet not the characteristic function of a set. This allows the IMF to retain some non-stationarity features.

\bigskip

\section{Application to time dependent signals}\label{sec:applTimeDep}

Consider a signal $\sigma$ on the time interval $[\,0,L\,]$ sampled at rate $B$ per time unit. This signal is represented as a vector of size $N=BL$. As it is customary in the signal processing literature, we wrap this signal on a discrete circle: consider the quotient group $\Z_N=\Z/N\Z$ and set $t_i=i/B$, $i\in \Z_N$. Our discretized and wrapped signal is $s=(\sigma(t_i))_{i\in\Z_N}$.

The discrete Fourier transform $\mathcal Fs$ is calculated on the points $k/L$, $k\in \Z_N$, in the interval $[\,0,B\,]$. Let $\omega$  be  $e^{-2\pi  \imath/N}$, the $n$-th root of unity. The Fourier transform operator $\mathcal F$ is the $N\times N$ matrix with entries ${\mathcal F}_{i,j}= \omega^{i-j}/\sqrt{N}$, with $i,j = 0, 1, \dots, N-1$.

Next we discuss three key ingredients of the FIF algorithm: the choice of the filters $w_j$, the powers $p_j$, and how to ensure that the algorithm stops in finite time.  We only give a brief discussion that suffices to convey the essential ideas, and refer the reader to \cite{cicone2021numerical, cicone2019MFIF} for details.

\subsection {Choosing the filters}\label{sec:choosingfilters}
As indicated at the end of last section,  in view of (\ref{eq:keyEq}),  we are more interested in the shape of $\mathcal{F}w$, the Fourier transform of the filter than the filter itself. Since $w$ is compactly supported, its Fourier transform has unbounded support.  Definition \ref{def: filter} (ii) implies that $\mathcal{F}{w}(0)=1$. For the particular filter that we pick, $\mathcal{F}{w}$  decreases from $\mathcal{F}{w}(0)$ until it reaches zero. After this smallest positive zero $\xi$,  it oscillates with smaller and smaller amplitude.
By requiring $\mathcal{F}w$ to be smaller than some desired quantity to the right of $\xi$ we create a damping effect in the frequency domain. By damping and not annihilating frequencies, the IMF retains the  non-stationarity of the signal, and yet concentrates the
frequency range of the IMF  above $\xi$.

We observe that the zero set of $\mathcal{F}w$ is preserved by self convolution of $w$ since $\mathcal{F}w^{\star k}=(\mathcal{F}w)^k$.  We use only double convolution filters.  This ensures that $\mathcal{F}w$ is in $[\,0,1\,]$ and $\mathcal{F}w^{\star k}$ decreases with $k$ except at the origin.
 Since $\ell(w^{\star k})=\min\{k\ell(w), N\}$, the filter length of $w^{\star k}$ increases with $k$. A long filter length makes moving averaging less local, and therefore increases stationarity.

In practice we choose a double convolution filter $w_0$ with support $[\,-1/2, 1/2\,]$
and take a dilation $w(\cdot)=w_0(\cdot /\ell)/\ell $.  To choose the  filter length and control the frequency focus,  we compute the median
of the distances between two consecutive local extrema. This median is an estimate of the half period of the  highest observed frequency. We take the filter length $\ell$ to be $\nu$ times this median, where $\nu$ is a tuning parameter that accounts for the concentration of the filter.

We follow the iterative construction of the IMFs described in the previous section. The $j$-th filter $w_j$ is constructed as a dilation of $w_0$ as indicated above from the local extrema of  the signal $s_{j-1}$.
Each filter focuses on sections of the spectrum of the signal.
The first one, $w_1$,  focuses on a group of highest frequencies, the second one, $w_2$, the next group of frequencies, etc.

If $\xi_{j}$ is the first positive zero for the $j$-th filter, the first $\text{IMF}$ focuses on frequencies around and higher than $\xi_{1}$, the second  $\text{IMF}$  focuses on frequencies around $\xi_{2}$, etc. Because of the damping effect, there is a small residue of frequencies that are higher than $\xi_{j}$ in the $j$-th IMF.  Numerically we have observed from the periodogram of the $j$-th IMF  that its energy  is concentrated around the frequency $\xi_{j}$, and very little beyond $\xi_{j-1}$. Put differently, IMFs tend to keep together nearby frequencies on the important parts of the periodogram.

 \subsection{Deciding the power $\bm{p_j}$  }\label{sec:decidingPj}

    Equation (\ref{eq:keyEq}) gives us an estimate for the norm of the matrix $\| \mathcal{V}_{w} ^{p+1} - \mathcal{V}_{w}^{p} \|$. Indeed, write $D$ for $\diag(1 - \mathcal{F}w)$.  The function $(1-x)^px$ being bounded by $1/(ep)$ on $[\,0,1\,]$, each diagonal element of $(1-D)^pD$ is between $0$ and $1/(ep)$. Since $\mathcal{F}$ is an isometry, we then have
 \begin{equation}\nonumber
  \| \mathcal{V}_{w}^{p+1} - \mathcal{V}_{w}^{p} \|= \|\mathcal{F}^{-1}(1-D)^pD\mathcal{F}\|\le \|(1-D)^{p}D\|\le 1/(ep).
\end{equation}
 In practice,  we set a threshold $\delta $, typically $10^{-3}$ or $10^{-4}$,  and for a signal $s$, we choose $p$ to be the smallest integer such that
$\| \mathcal{I}_{w,p+1}s - \mathcal{I}_{w,p}  s \|_2\le \delta\|s\|_2$. The actual number $p$ is much smaller than that from the theoretical estimate.

  \subsection{Ensuring the algorithm terminates}

  Because IMFs are not bandpass filters but instead just focus on sections of the spectrum,  the remainder signal $s_{j}$ may still contains some high frequencies. This makes possible for the
  filter length chosen as in section \ref{sec:choosingfilters} to decrease. To prevent such a decrease we  force the filter length to increase by $10\%$ or some appropriate amount of the previous filter length when this happens.

  It is also possible that an IMF is insignificant, that is,   $\|\mathcal{I}_{w, 1}  s_{j}\|_2\le \delta\|s_{j}\|_2$.  In this case, we redo the iteration by forcing the filter length to increase by $10\%$. Since the signal is on the discrete circle, its frequency domain contains finitely many points. This $10\%$
  increase guarantees that the algorithm stops. In real life applications that we encountered, there are sufficiently many significant IMFs, and only a trend with at most one extrema is left as the last remainder.

\section{IMFogram}

The periodogram and its localized version, the spectrogram, are plots adapted to Fourier decomposition. Because each Fourier component has a specific frequency and these components are mutually orthogonal, the total energy of the signal is the sum of each component's energy. In contrast, each IMF does not have a fix frequency and the IMFs are not mutually orthogonal. The IMFogram is a simple analog of the spectrogram that can be quickly computed on IMF decompositions. It is defined as follows:

Let $\eta$ be a parameter --- $5$ to $20
$ in practice.
Consider an IMF, denoted by $f$, which was produced with a filter of length $\ell$.  Its local energy at $t$ is approximated by
 \begin{equation}\nonumber
  E_f(t) =\frac{1}{2\eta \ell}\int_{t-\eta\ell}^{t+\eta\ell} f(\tau)^2 \d\tau.
\end{equation}

Because an IMF focuses on a narrow frequency section, the local frequency of $f$ may be approximated by
 \begin{equation}\nonumber
 \Omega_f(t) = \frac{1}{4\eta\ell} \times \text{number of $0$-crossings of $f$ over the interval  $[\,t-\eta\ell, t+\eta\ell\,]$}.
\end{equation}
There are many other possible approximations discussed comprehensively in \cite{huang2009convergence}.  We obtain time, frequency, energy triples $(t_i, \Omega_{f}(t_i), E_{f}(t_i))_{i\in\mathbb{Z}_N}$ for $f$.

Consider a discretized signal $s$ as in section \ref{sec:applTimeDep}. Its time domain in $\{i/B: i\in\mathbb{Z}_N\}$ and its frequency domain is $(\mathbb{Z}_N/L)\cap[\,0,B/2\,]$. We partition the time-frequency domain in rectangles. For a rectangle $R$, let $\Pi_{t}R$ be its projection onto the time coordinate and $\Pi_{\omega}R$ its projection onto the frequency coordinate. Since the time domain is discrete, the cardinality $\# \Pi_t R$ is finite.
Let $s_1, s_2, \dots, s_k$ be the IMF decomposition of $s$. To a rectangle $R$ we associate an energy defined as a sum of the average local energies of each IMF when the local frequency lies in $\Pi_{\omega}R$, that is
 \begin{equation}\nonumber
 E_s(R) = \sum_{1\le j \le k} \frac{1}{\#\Pi_{t}R}\sum_{\tau\in\Pi_{t}R}  E_{s_j}(\tau)\mathds{1}\Bigl\{\Omega_f(\tau)  \in \Pi_{\omega}R \Bigr\}
\end{equation}
We recommend to choose the length of $\Pi_{t}R$ comparable to the smallest filter length $\ell(w_1)$ so that changes in high frequencies are well represented. Depending on the application, one may take rectangles of different sizes to cover the time and frequency domain.

The IMFogram of $s$ is the plot of the step function that equals to $E_s(R)$ on each rectangle $R$. Sometimes one may want to discard small IMFs when plotting the IMFogram for a signal.

\section{Numerical Examples}\label{sec:NumericalExamples}

In this section we present a few examples of time frequency analysis of both synthetic and real world signals, and the ability of FIF to handle higher dimensional data. The double convolution filter $w_0$ that we use for these examples is obtained by taking $\tilde{w}$ to be the stationary solution from a  Fokker-Plank equation \cite{cicone2016adaptive}. This approach shares similarities with \cite{wang2012iterative} where filters are solutions of a Fokker-Plank equation, and filter length is in one-to-one correspondence with the time in the evolution equation. Our approach results in a faster algorithm because the filters are dilation of a fixed $w_0$.

\subsection{Synthetic signal}

\begin{figure}
	\centering
	\includegraphics[width=0.5\linewidth]{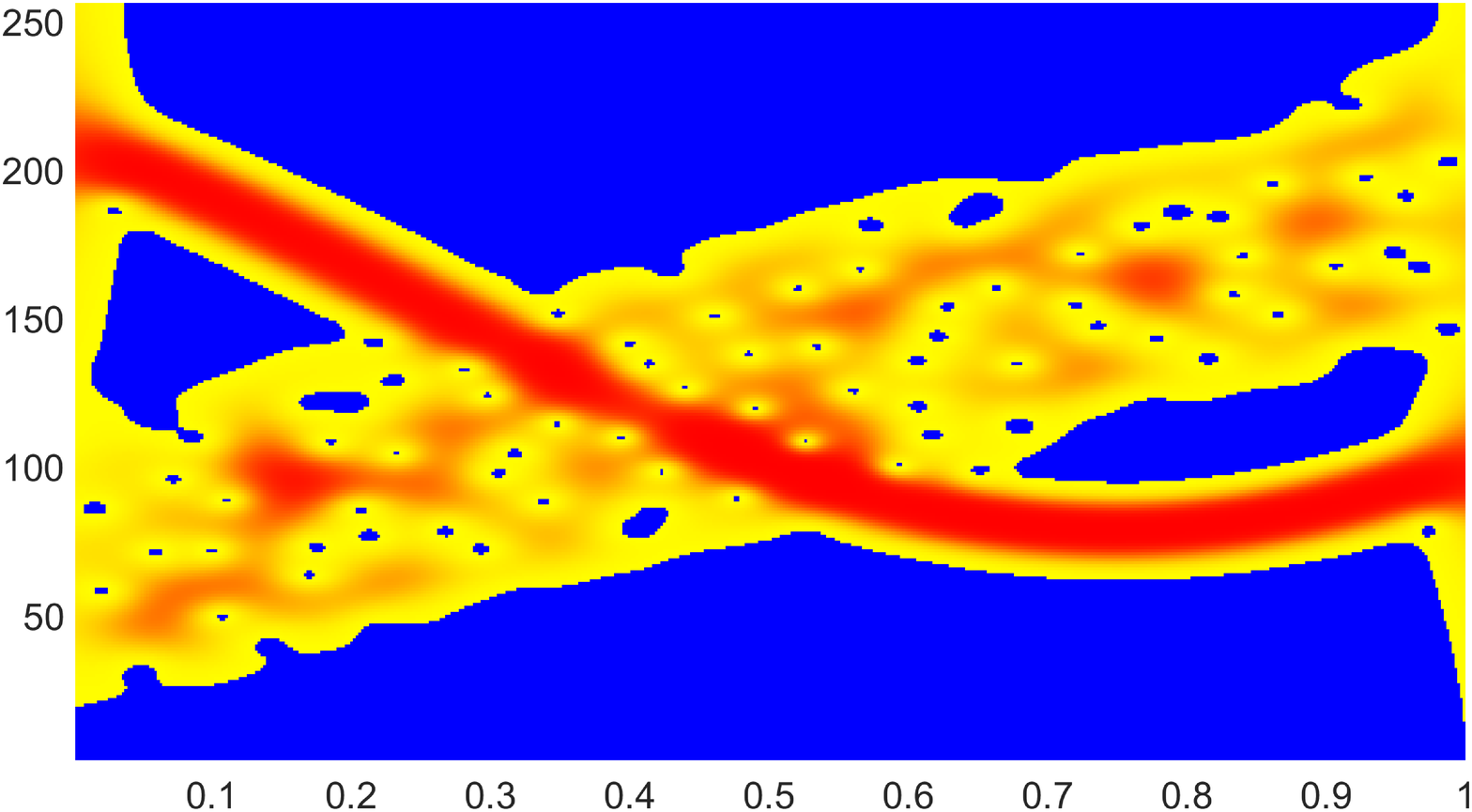}~\includegraphics[width=0.5\linewidth]{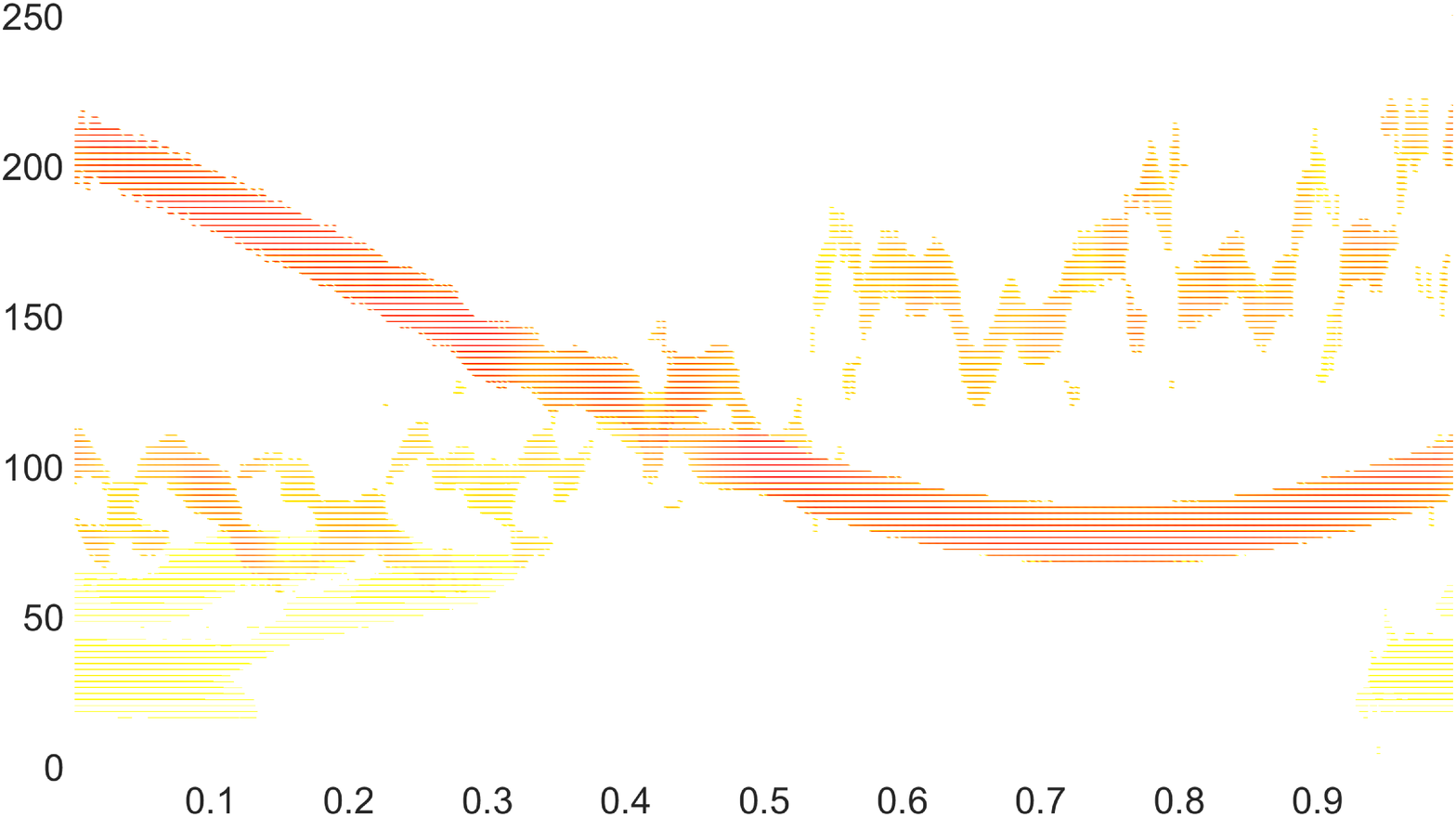}\\
	\includegraphics[width=0.5\linewidth]{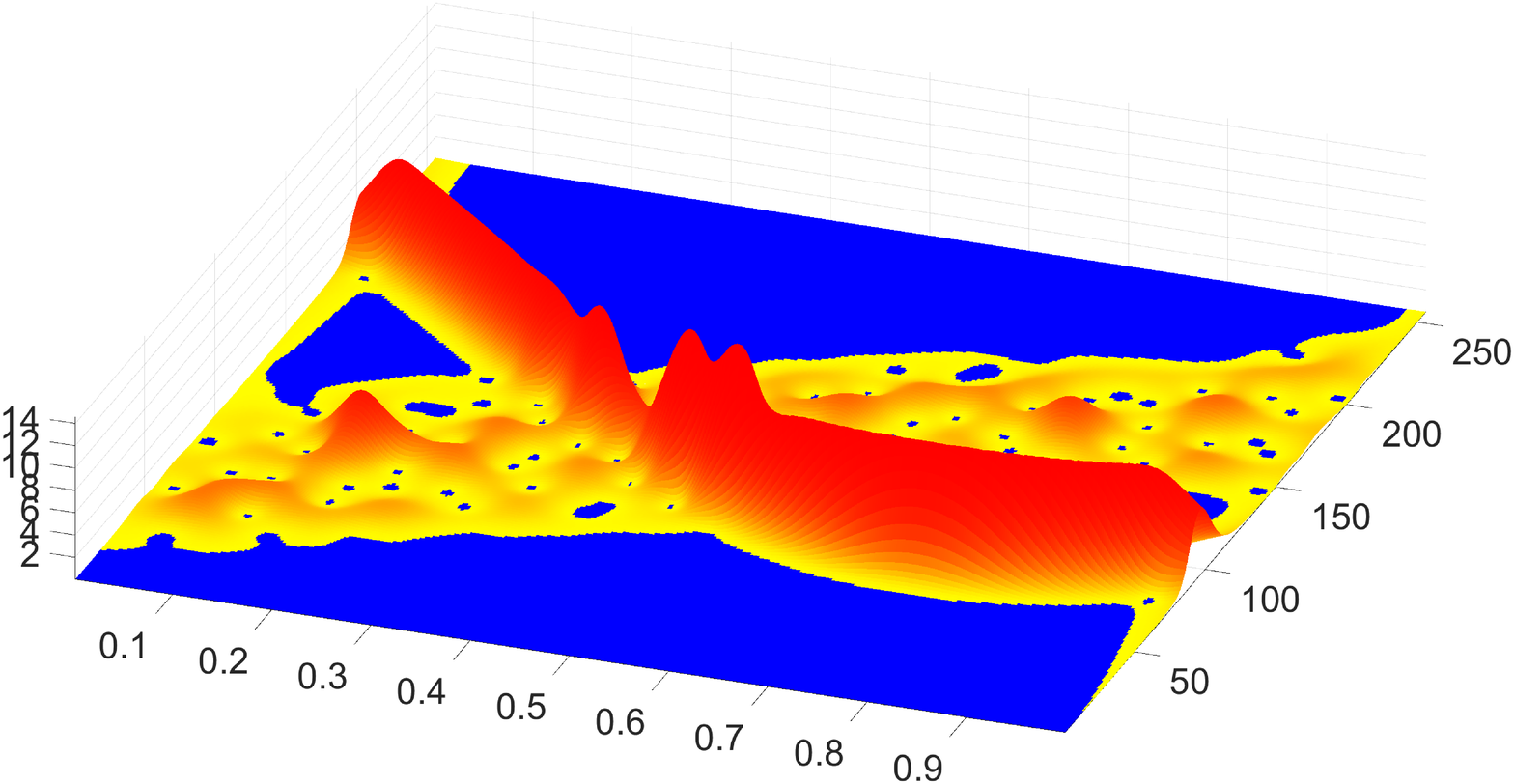}~\includegraphics[width=0.5\linewidth]{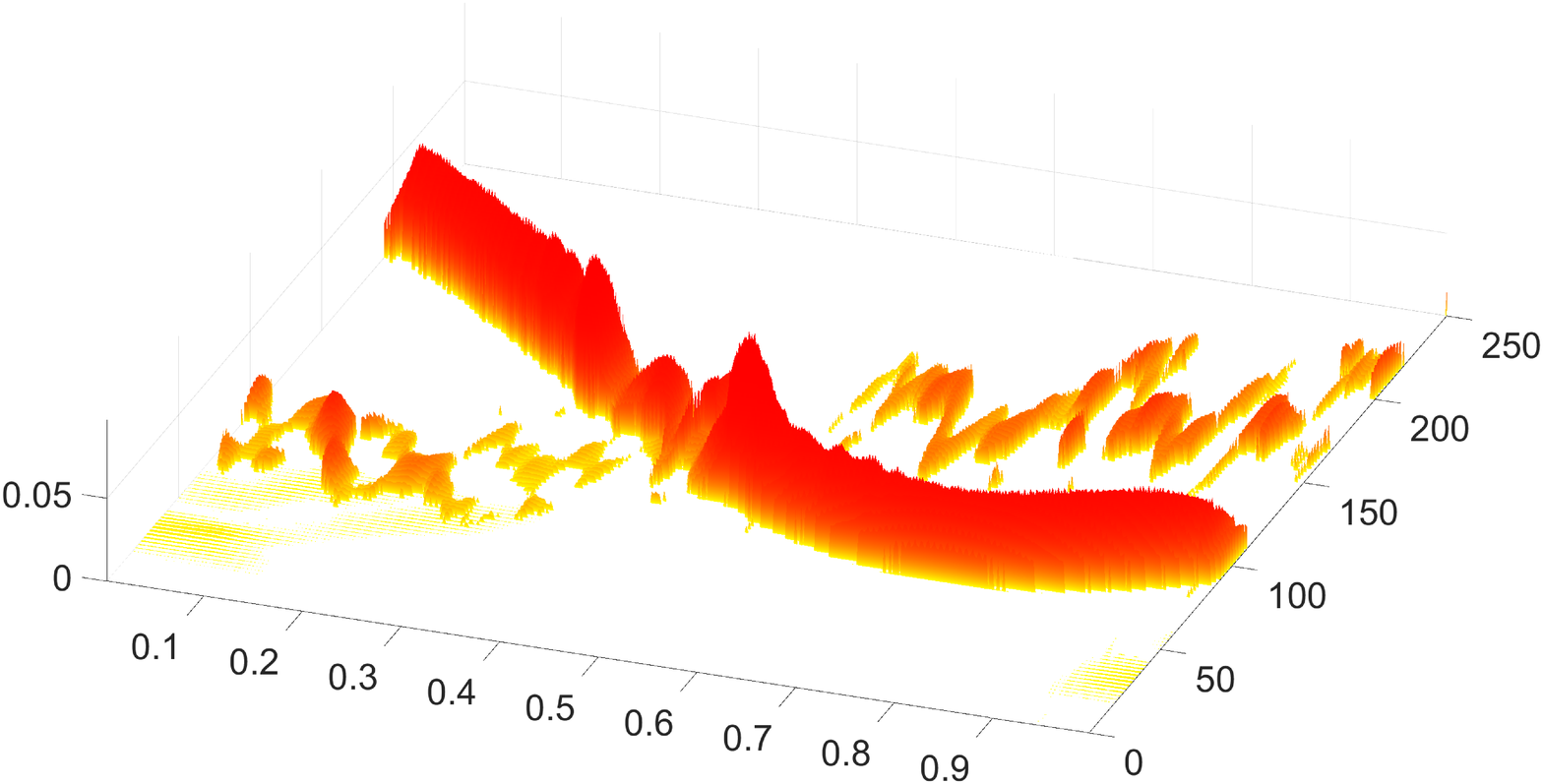}
	\caption{Time frequency representation of one chirp in time-varying noise. Left: spectrogram. Right: IMFogram.}
	\label{fig:Ex1_TFA}
\end{figure}

Consider a synthetic signal mixing a nonlinear chirp (rapid change in frequency)  and a bandpass time-varying noise, Figure \ref{fig:Ex1_TFA}. This example was studied in other previous works, like for instance in \cite{xiao2007multitaper}.

 Figure \ref{fig:Ex1_TFA} compares the classical short time window spectrogram and IMFogram. As we see, the IMFogram is more concentrated on the time-frequencies used to generate the signal.

\subsection{Piano recording}

  \hskip .80 \textwidth \raisebox{-1.8 cm}[0pt][0pt]{%
    \includegraphics[width=3cm,height=4.5\baselineskip]{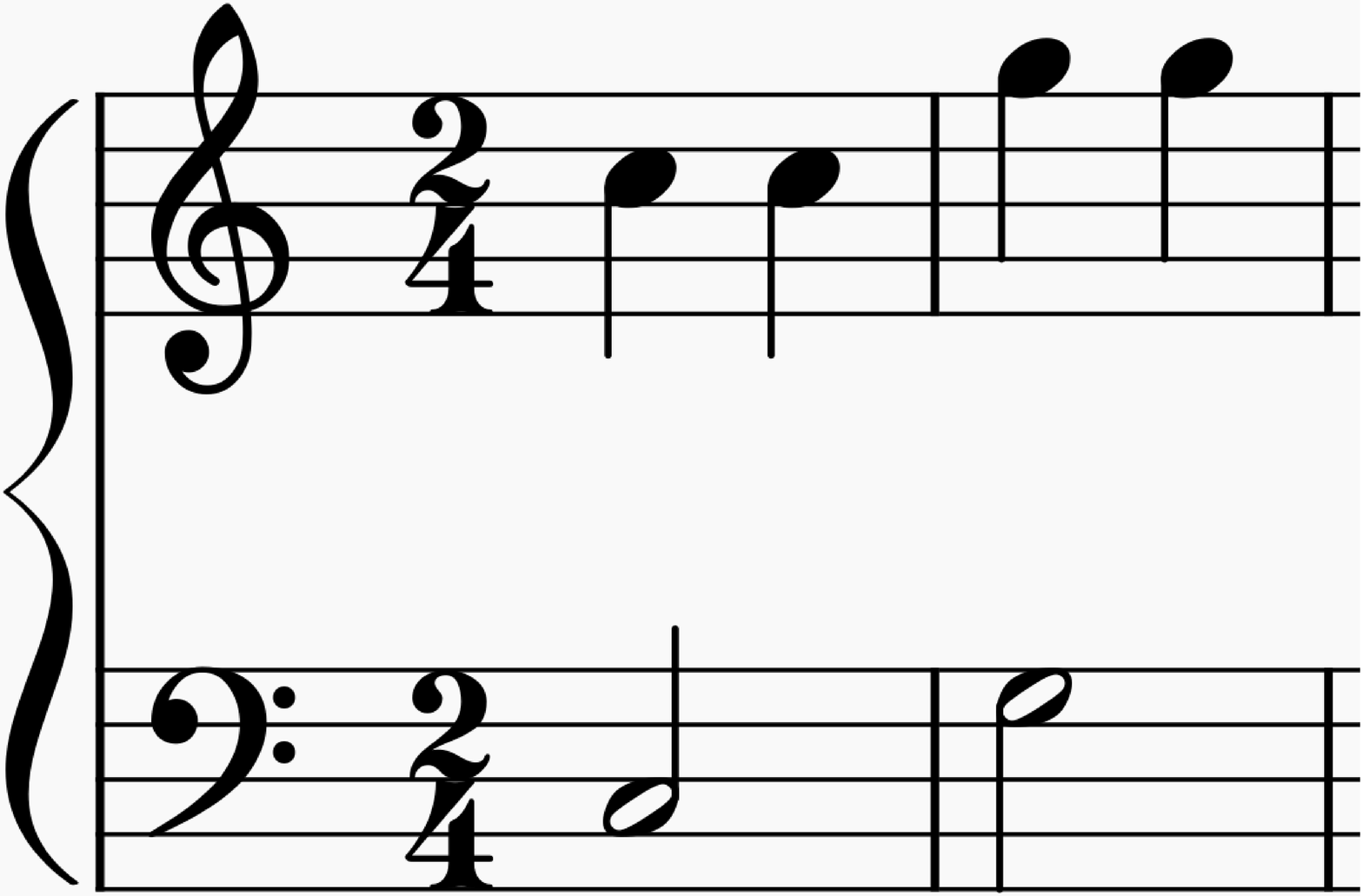}}
\vskip - \baselineskip
\parshape 2
  0pt 0.75\textwidth
  0pt 0.75\textwidth
\noindent We consider the time frequency analysis of the first four notes of {\it A vous dirais-je Maman}.

\begin{figure}
	\centering
	\includegraphics[width=0.5\linewidth]{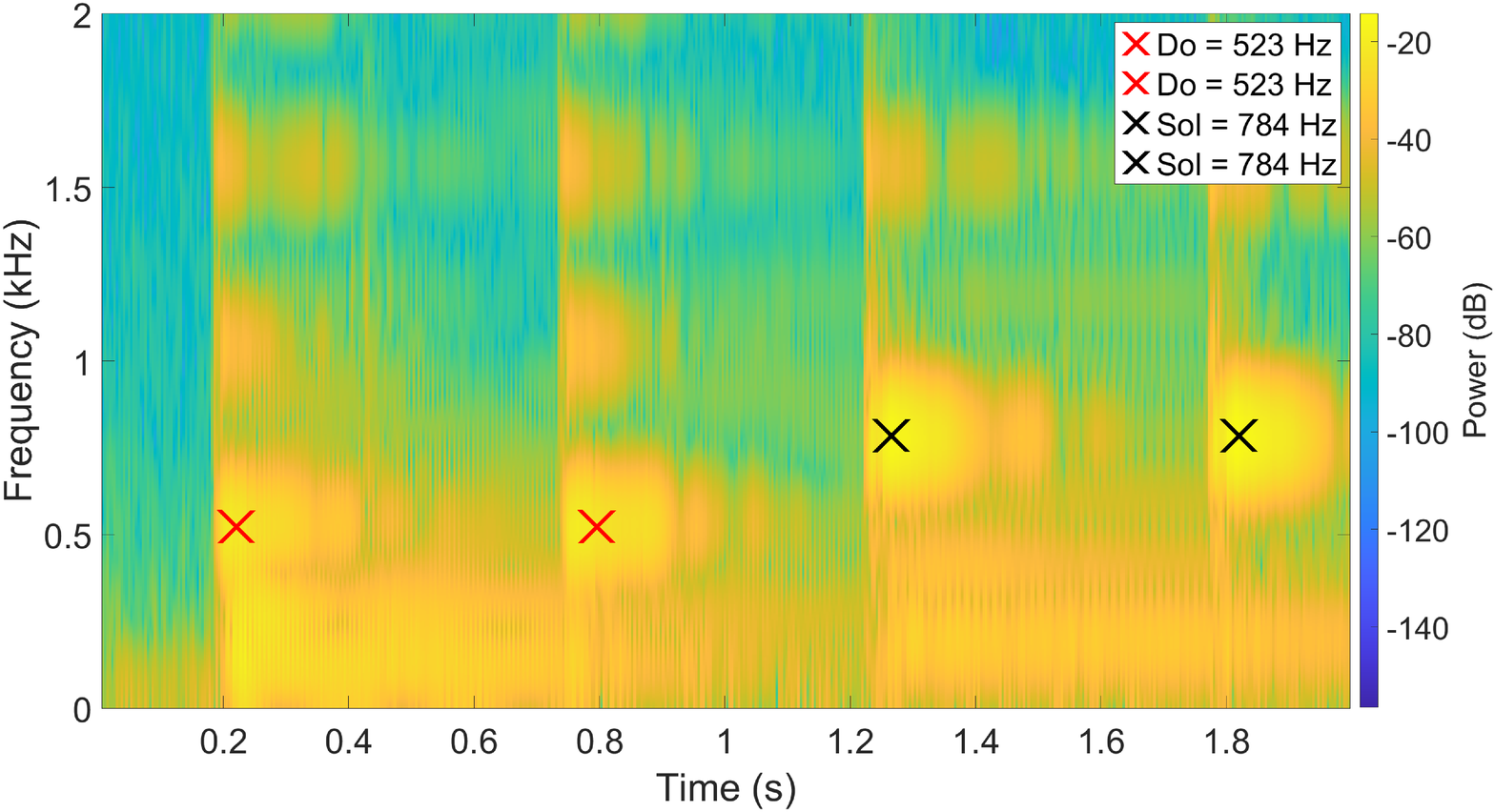}~\includegraphics[width=0.5\linewidth]{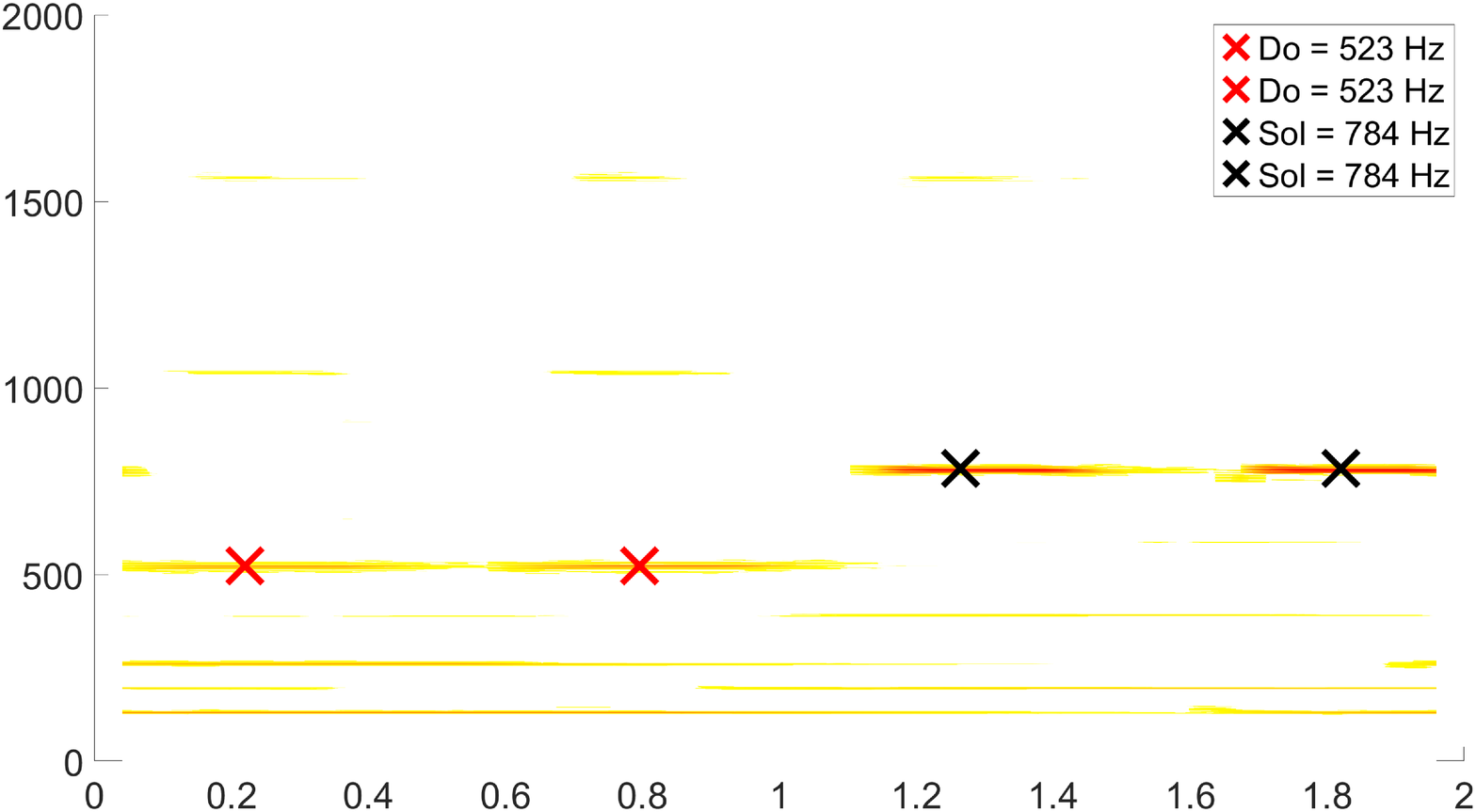}\\
	\includegraphics[width=0.5\linewidth]{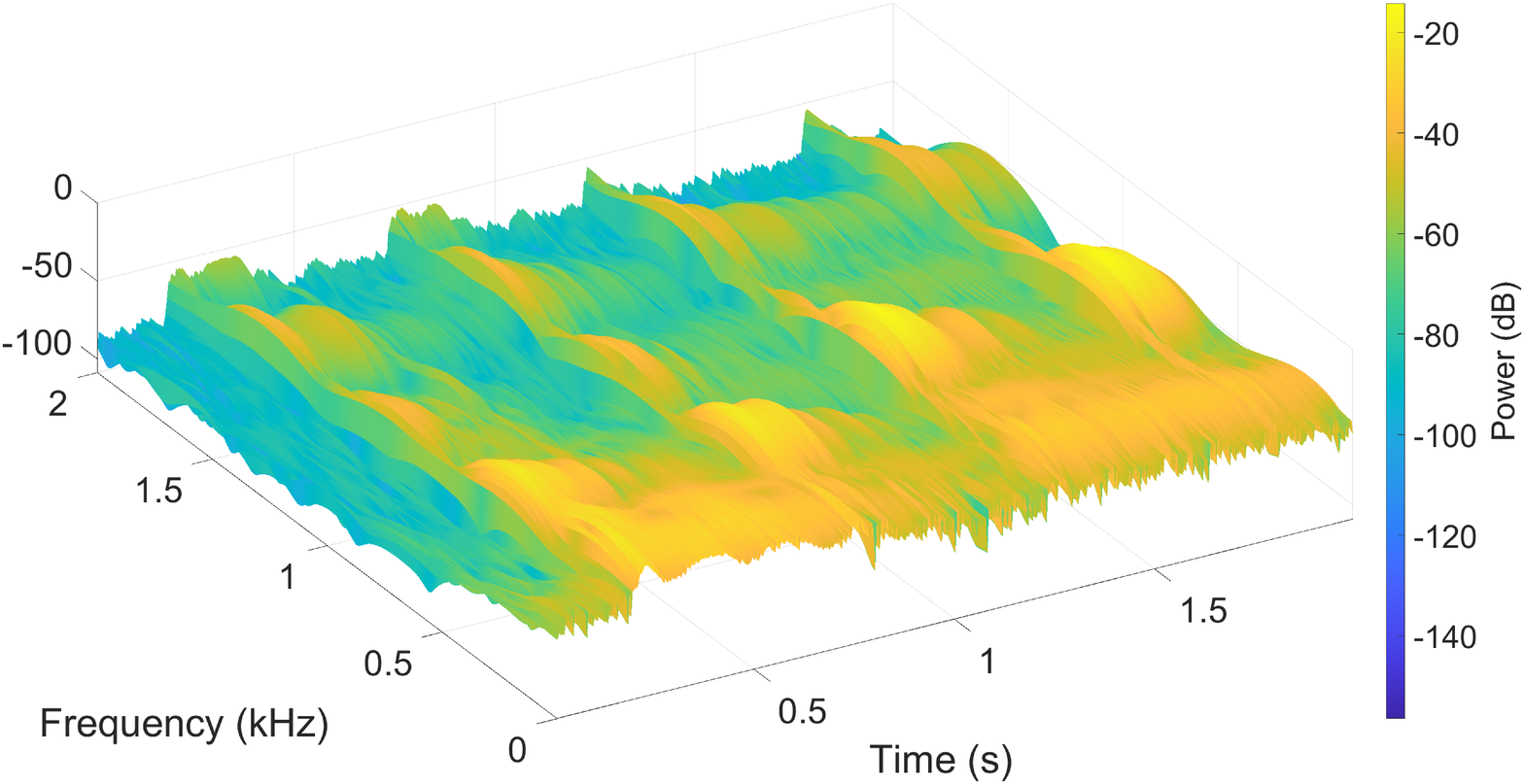}~\includegraphics[width=0.5\linewidth]{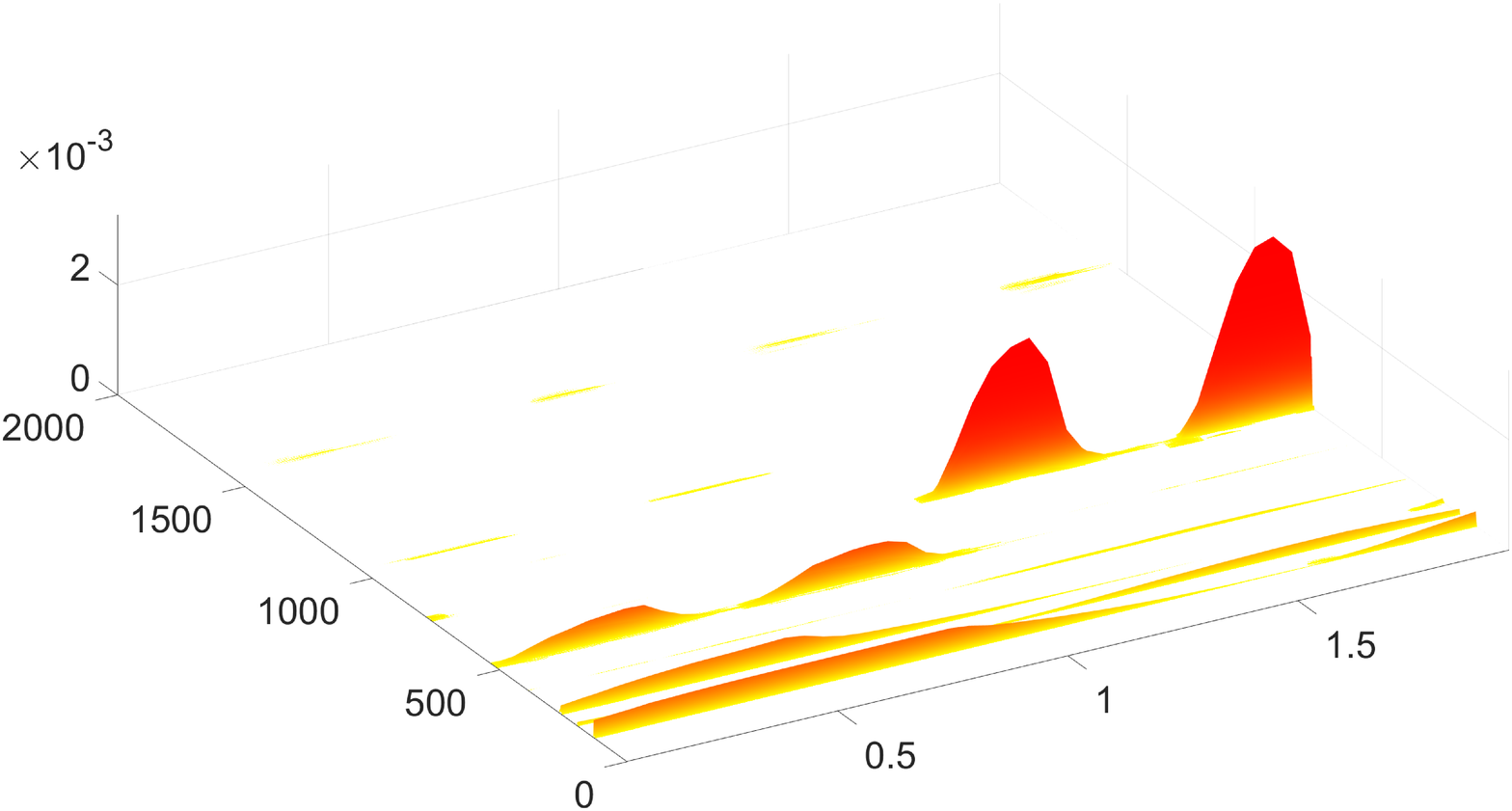}
	\caption{Time frequency representation of the first four quarter notes of the {\it A vous dirais-je Maman}.   Left: Spectrogram. Right: IMFogram.}
	\label{fig:Ex2_TFA}
\end{figure}

\begin{figure}
	\centering
	\includegraphics[width=0.71\linewidth]{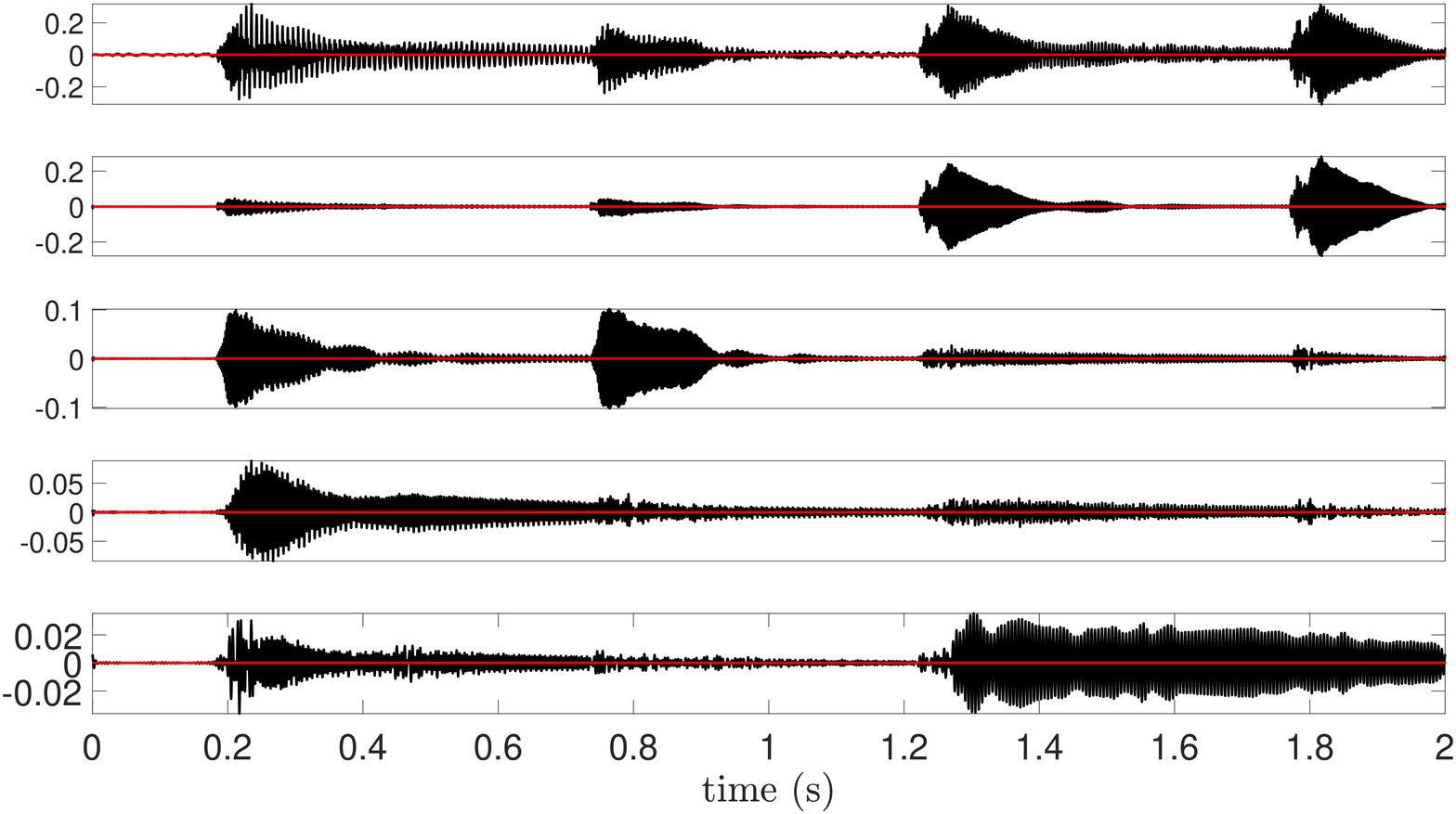}~\includegraphics[width=0.28\linewidth]{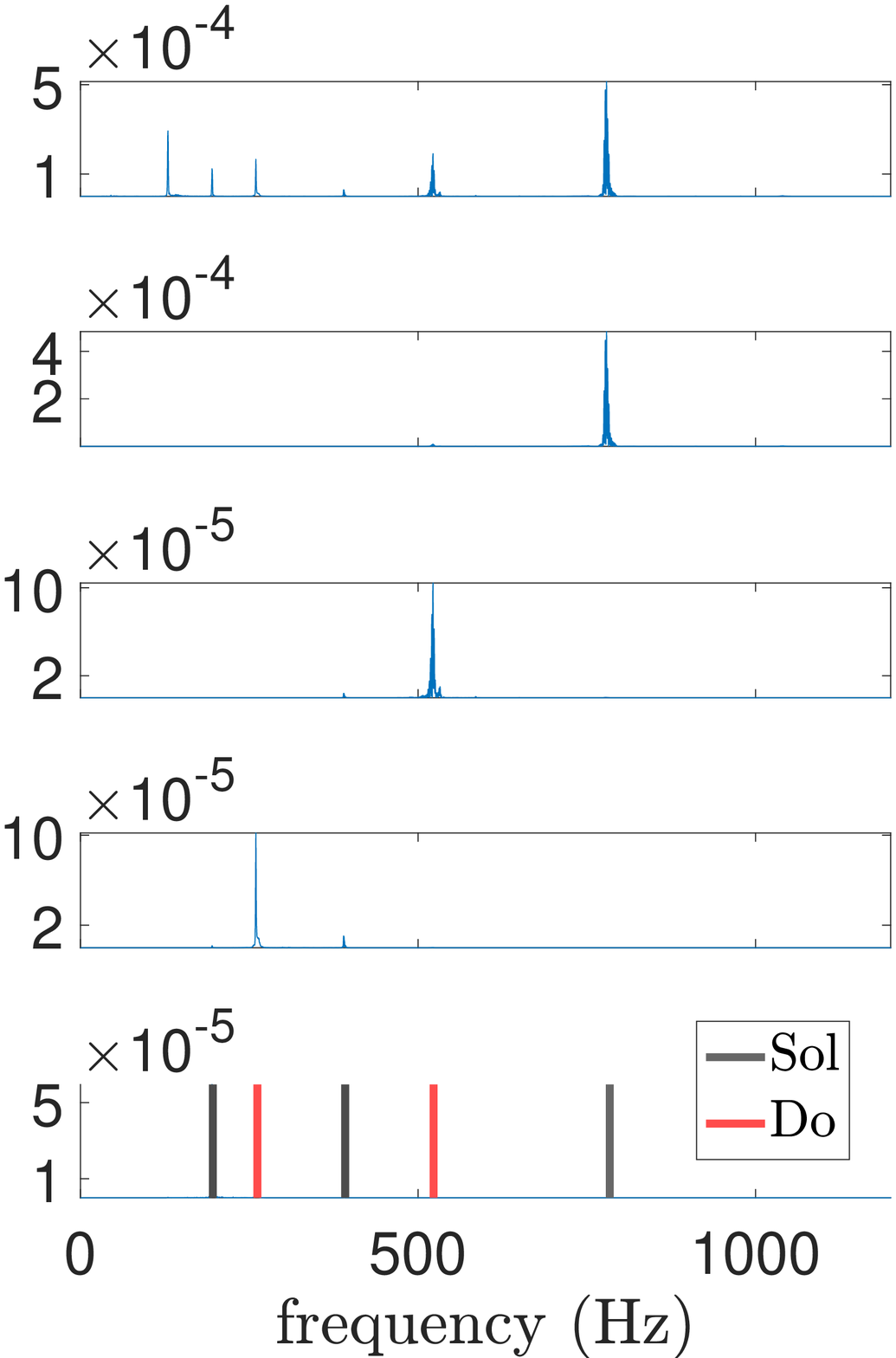}
	\caption{Top row: original piano recording and its periodogram.  Next four rows: 4-th to 7-th IMFs and their periodograms.}
	\label{fig:Ex2_FIF_decomp}
\end{figure}

\parshape 3
  0pt 0.75\textwidth
  0pt 0.75\textwidth
  0pt \textwidth
The four notes from the right hand are Do Do Sol Sol in quarter notes and the left hand are Do Sol in half notes. Figure \ref{fig:Ex2_FIF_decomp} shows the original recording, top row, and the $4$-th, $5$-th, $6$-th, and $7$-th IMFs, and their periodograms. We see that the $4$-th IMF pick up the upper Sol, the next one the Do beat opening the piece. The left hand is captured by the $6$-th and $7$-th IMFs. We also see some harmonics of the notes, which reflects the fact that piano notes are complex sounds, especially on the lower register. Compare with spectrogram, we see in  Figure \ref{fig:Ex2_TFA} that the IMFogram gives a sharper view of the notes. This example suggests that FIF may be of value for automatic scoring, which is useful for music transcription and its applications.

\subsection{Global temperature}

\begin{figure}
	\centering
	\includegraphics[width=0.5\linewidth]{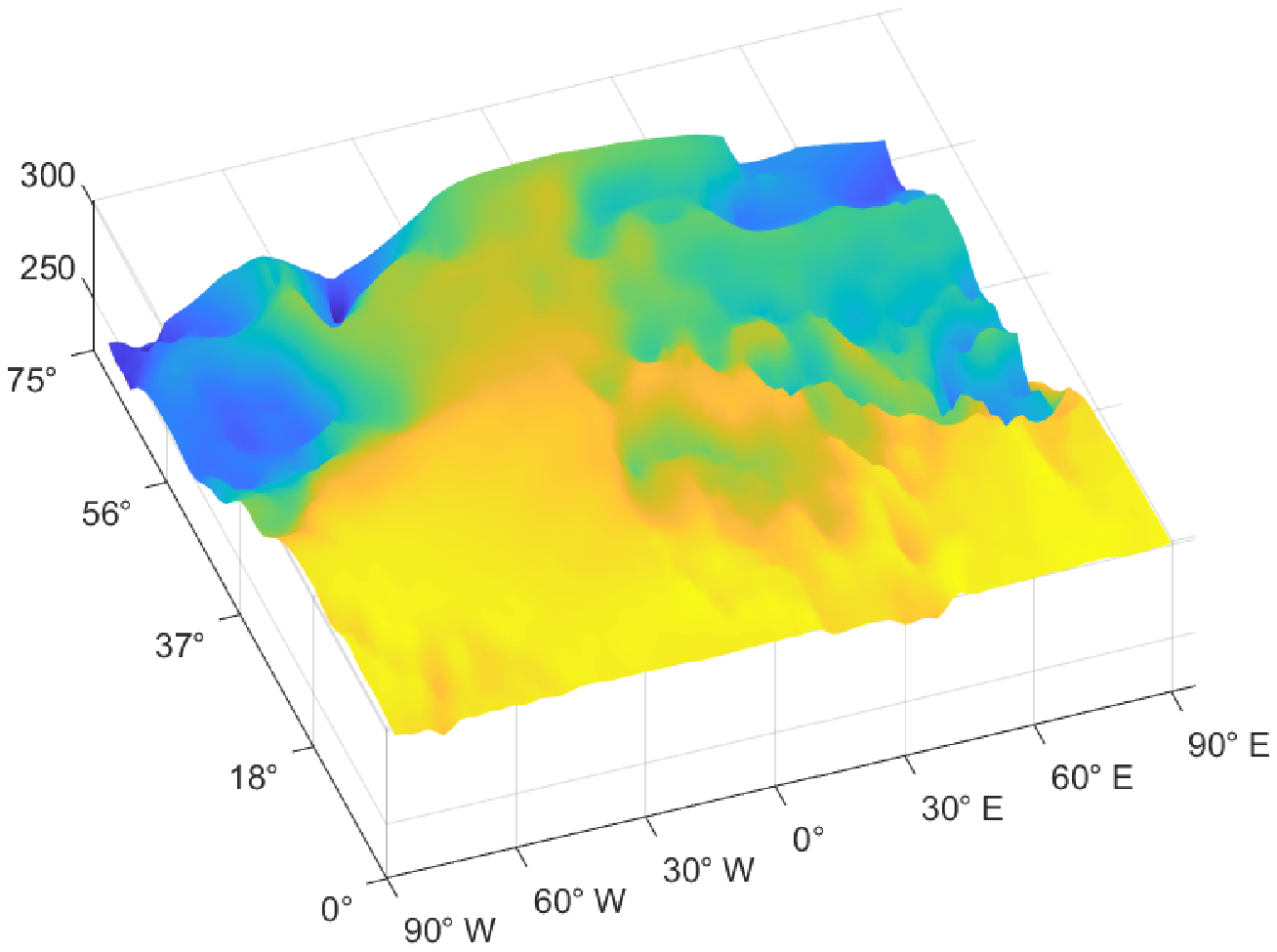}~\includegraphics[width=0.5\linewidth]{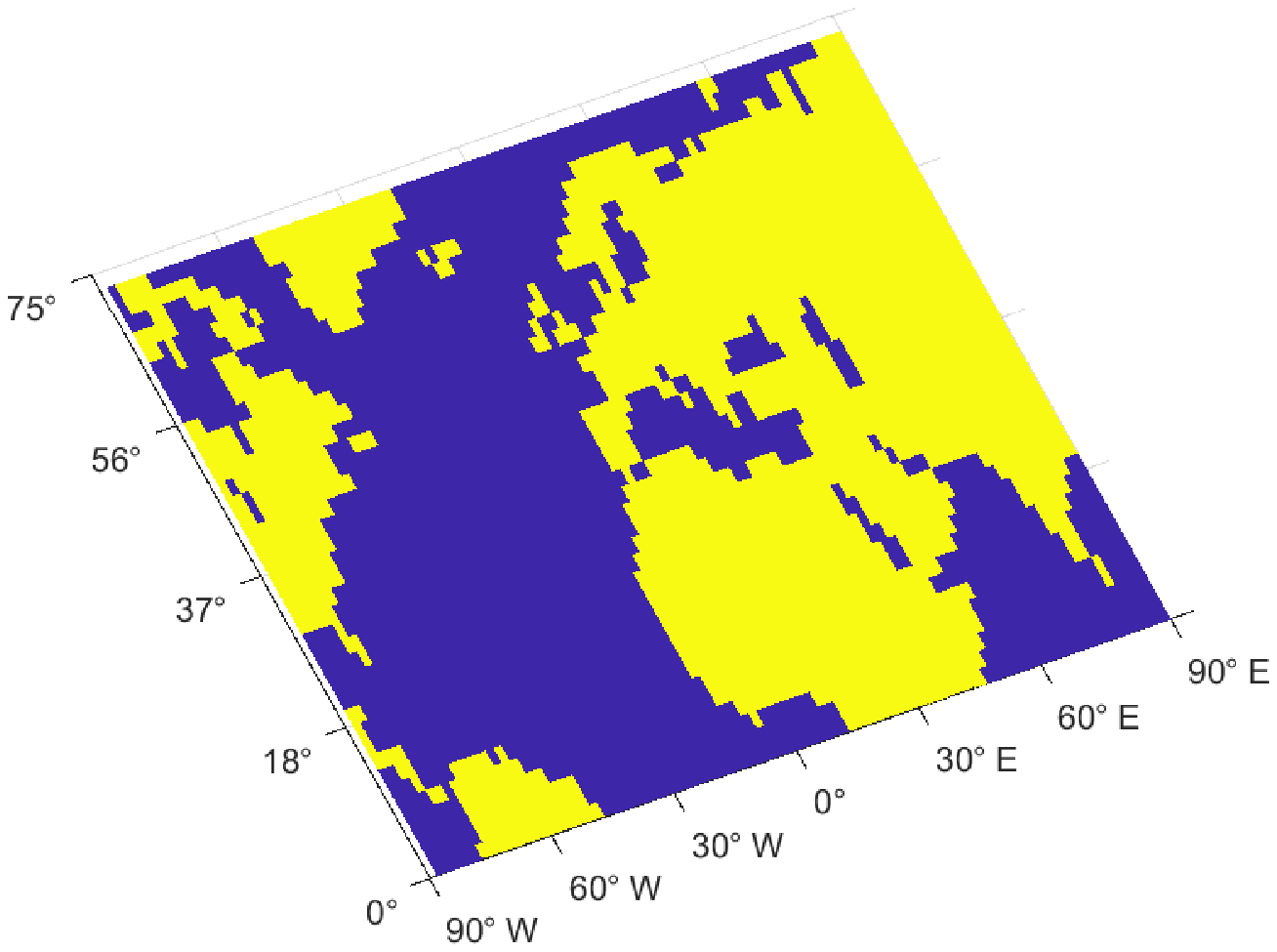}\\
	\includegraphics[width=0.5\linewidth]{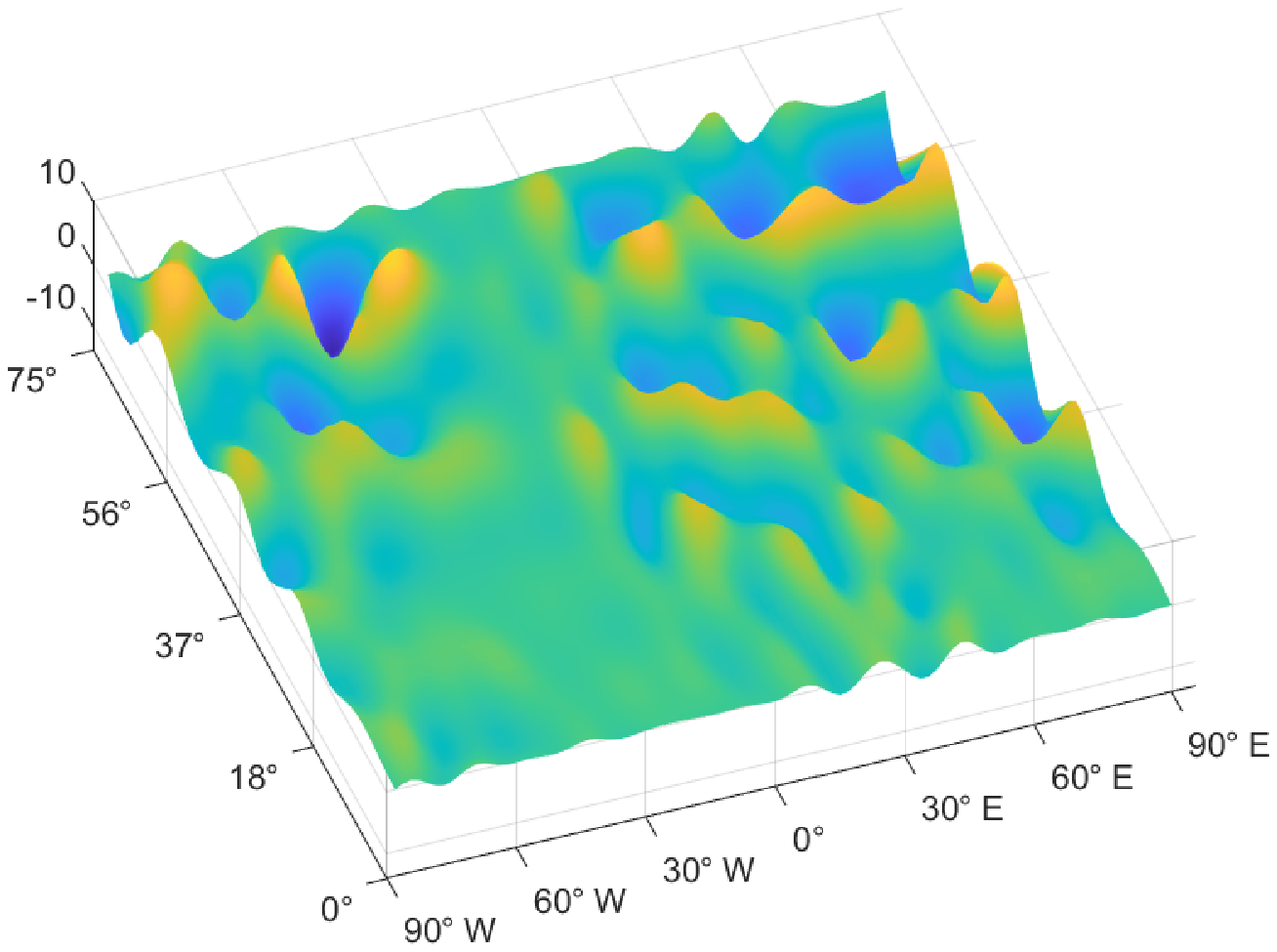}~\includegraphics[width=0.5\linewidth]{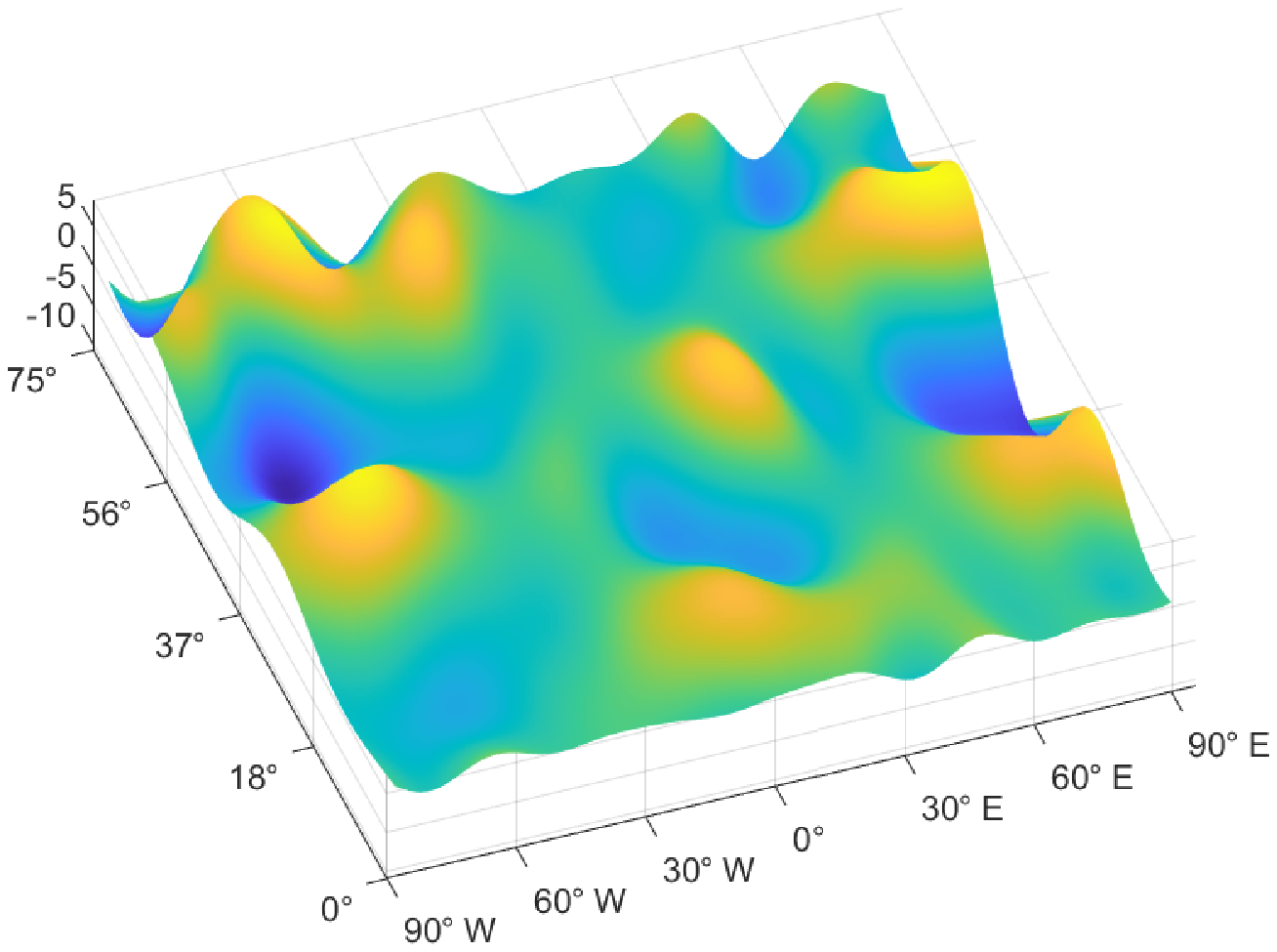}\\
	\includegraphics[width=0.5\linewidth]{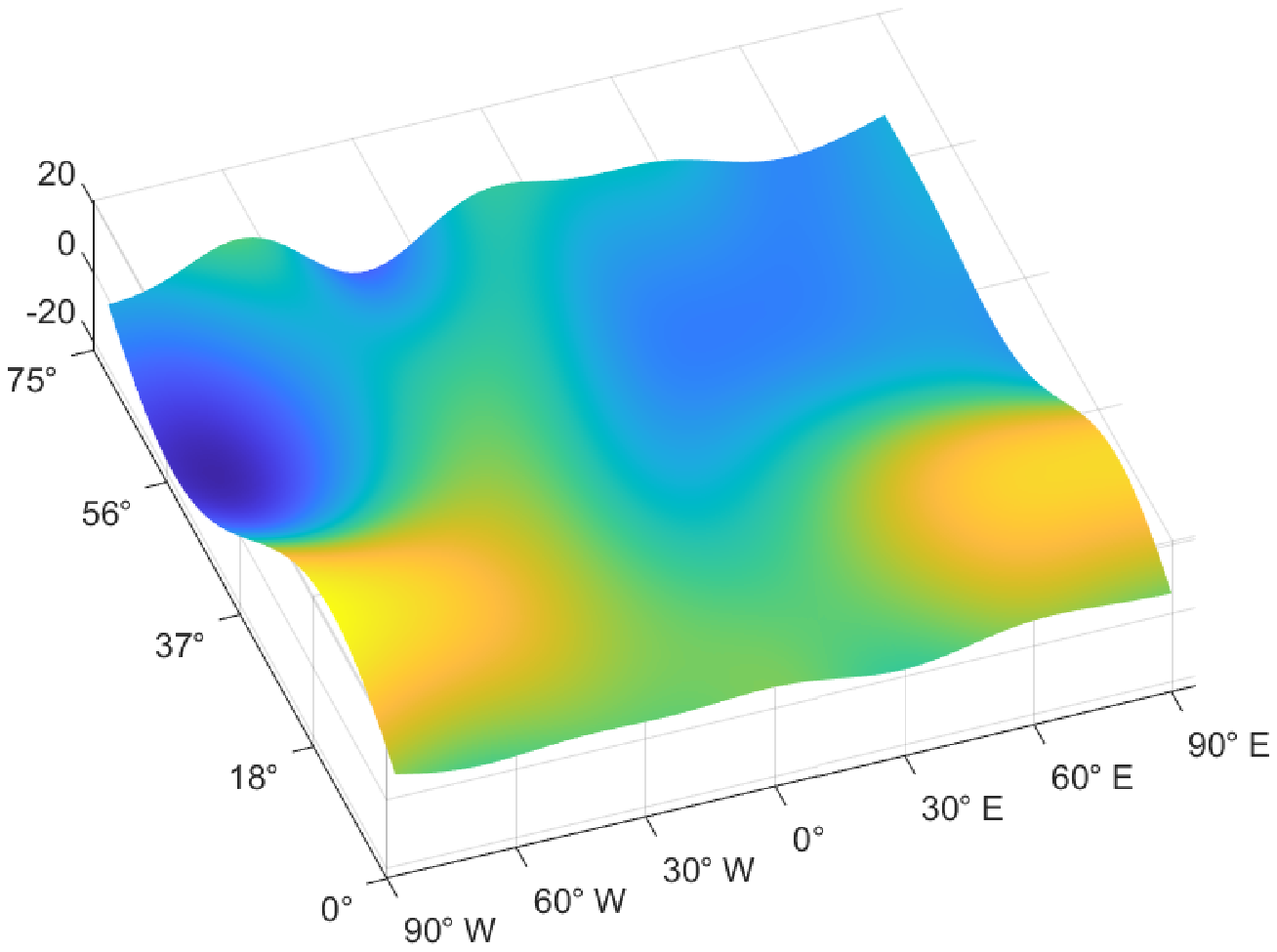}~\includegraphics[width=0.5\linewidth]{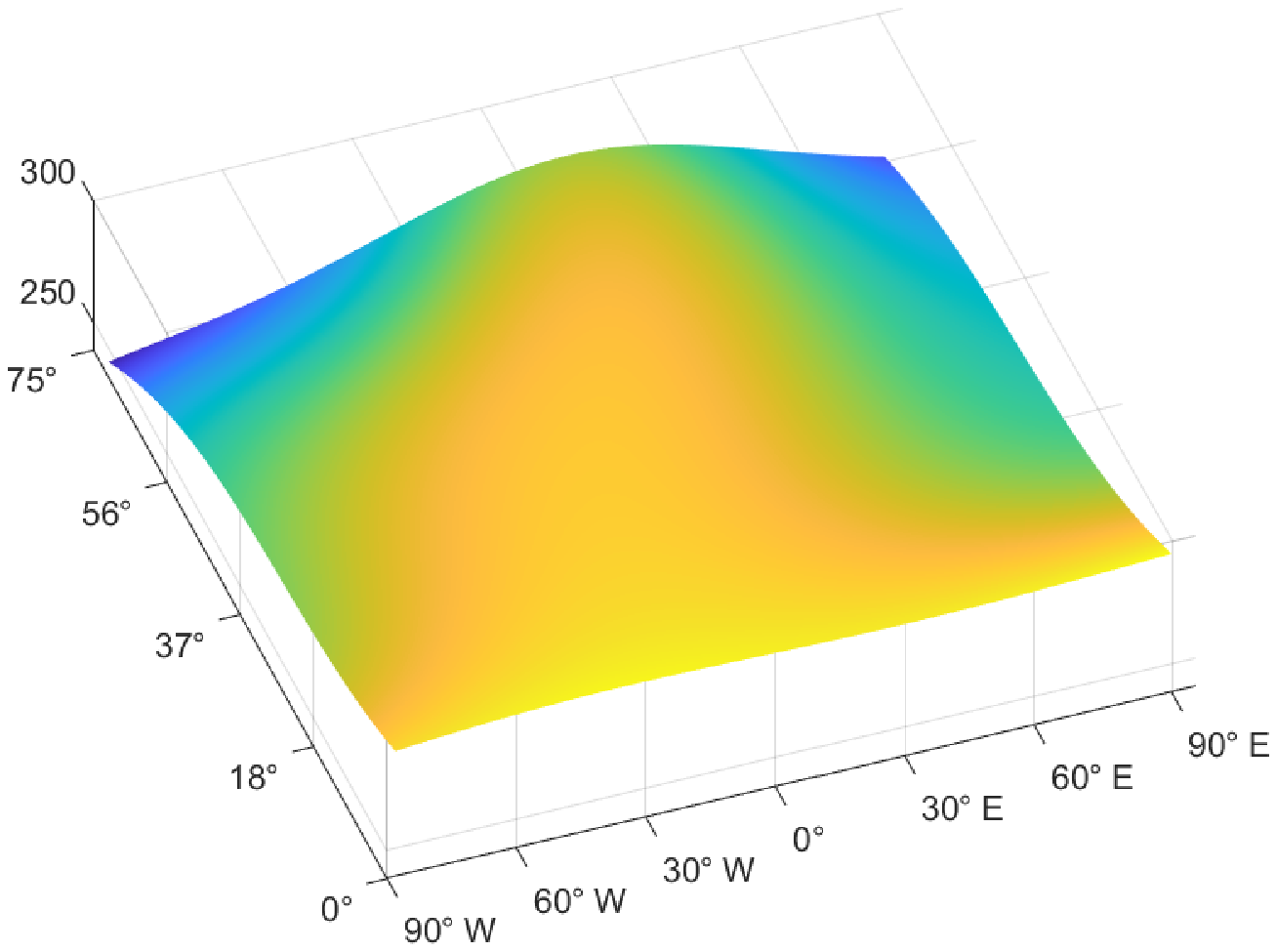}
	\caption{Top left: Earth's  air temperature (in K) on January 1, 2014, in a quarter-spherical window covering Eastern America and West/Central Eurasia. Top right: land and sea map corresponding to the data. Middle and bottom rows: some of the IMFs components generated by the FIF2 algorithm.}
	\label{fig:Ex3}
\end{figure}

To show the ability of FIF to handle $2$-dimensional data. The Earth's air temperature measurements are made available through the NCEP/NCAR Reanalysis project. The original NCEP Reanalysis data have been provided by the NOAA/OAR/ESRL PSD, Boulder, Colorado, USA  \cite{NOAA} and have been edited by the Climatic Research Unit, University of East Anglia \cite{UEA}. They consists of measurements taken on a global grid at an altitude of 2 meters from the surface in Kelvin.
We choose the data from January 1, 2014, for a quarter-spherical window between $0^{\circ} $N -- $75^{\circ} $N and $90^{\circ}$W -- $90^{\circ} $E. This covers eastern North America, most of the north Atlantic Ocean, Europe, Africa above equator, western and continental Asia, as illustrated in the top right panel of
Figure \ref{fig:Ex3}.

Some of the IMFs produced by a two dimensional FIF method \cite{sfarra2021maximizing,cicone2017multidimensional} are exhibited on the second and third rows of Figure \ref{fig:Ex3} show.
The bottom right panel shows a broad trend: temperatures are about constant and warm on the equator, decreasing on the more continental land masses of the northern hemisphere, with the Atlantic ocean warmer than those northern land masses.
The bottom left panel shows broad fluctuations around this trend. The continental climates of northern USA and Canada, similar to that of Russia, are markedly colder than other regions.
On the middle right panel, the depression near $90^{\circ}$E --  $31^{\circ}$N corresponds to the Himalaya region, while that around $85^{\circ}$W -- $50^{\circ}$E corresponds to the north of the Great Lakes area.
The middle left panel gives other fluctuations that may have other interpretations.


\section{Conclusion}

The examples show the value of the IMFogram and confirm that of iterative filtering methods.

The IMFogram is developed as a visualization tool to ``see'' the local frequency and amplitude information of a signal simultaneously. By varying the partition of frequency-amplitude space, it provides the visualization in multiple resolutions.  Not only the IMFogram can be used together with the FIF and EMD algorithms as demonstrated in this paper, but also with other decomposition strategies for non-stationary signals.

By bringing sharper frequency capability than the spectrogram, the IMFogram seems a good candidate to improve automated frequency identification techniques.
The similarities and differences between spectrogram and IMFogram shown in the examples presented in this work suggest that there is a mathematical connection between these two time-frequency representations. We plan to explore this direction of research in a future work.

The group theoretic framework that we introduced shows that time indexed and spaced indexed methodologies can be unified and that iterative filtering methods can be applied to a wider range of problems than previously considered. This framework also raises interesting questions on discretization and sampling on homogeneous spaces and Lie groups which we hope to tackle in a future work.
	

\end{document}